\documentclass{article}
\usepackage{harvard,amsmath}
\usepackage{epsf}

\newtheorem{Theo}{Theorem}
\newtheorem{Coro}{Corollary}

\DeclareMathOperator{\Var}{Var}
\input{tcilatex}

\begin{document}

\title{Optimal Asset Allocation\\
with Asymptotic Criteria}
\author{Vladislav Kargin\thanks{%
Cornerstone Research, 599 Lexington Avenue, New York, NY 10022, USA;
slava@bu.edu}}
\maketitle

\begin{abstract}
Assume (1) asset returns follow a stochastic multi-factor process with
time-varying conditional expectations; (2) investments are linear functions
of factors. This paper calculates asymptotic joint moments of the logarithm
of investor's wealth and the factors. These formulas enable fast computation
of a wide range of investment criteria. The results are illustrated by a
numerical example that shows that the optimal portfolio rules are sensitive
to the specification of the investment criterion.
\end{abstract}

\section{\protect\bigskip Introduction}

\label{Introduction}

The asset returns are predictable: see \cite{cochrane99} and \cite%
{campbell_lo_mackinlay97} who overview the empirical support for asset
return predictability\footnote{%
For example, see empirical studies by \cite{balvers_cosimano_mcdonald90}, %
\cite{breen_glosten_jagannathan90}, \cite{campbell87}, \cite{cochrane91}, %
\cite{fama_french89}, \cite{pesaran_timmermann95}, \cite%
{pesaran_timmermann00}, \cite{cooper99}.}, and \cite%
{brennan_schwartz_lagnado97}, \cite{campbell_viceira99}, and \cite%
{bielecki_pliska99} who study optimal portfolio allocation in the situation
with predictable returns. Conventional dynamic stochastic programming needs
a lot of computing power to find optimal allocations, and a possible
approach to this problem is to choose an approximate investment criterion
that would make the problem numerically tractable. This paper studies
portfolio optimization for a class of such criteria and tests the
sensitivity of optimal portfolio rules to the choice of the investment
criterion.

In early 1970s Samuelson and Merton derived Hamilton-Jacobi-Bellman equation
for the solution of the dynamic asset allocation problem, an equation that
is valid for any type of return process and any type of investment criteria.%
\footnote{\cite{samuelson69}, \cite{merton69}, \cite{merton71}, and \cite%
{merton73a}} This partial differential equation is, however, hard to solve.
To circumvent this difficulty, \cite{bielecki_pliska99}\footnote{%
See also \cite{stettner99}, \cite{bielecki_pliska_sherris00}, \cite%
{bielecki_hernandez_pliska99}, \cite{bielecki_pliska01}} introduced the
risk-sensitive investment criterion, which makes the optimization problem
easier. The criterion assumes that the investor's preferences depend only on
the growth rate of the expected logarithm of portfolio value and on the
growth rate of its variance. This assumption appears to be excessively tight
since it ignores other relevant moments, such as the covariance with the
factors. An investor might be concerned with the performance of the
portfolio in ``bad'' times, when his other sources of income are low. The
criterion that depends only on variance and expectation of the portfolio
value ignores this consideration.

This paper uses a new method to compute the asymptotic joint moments of the
logarithm of the portfolio value and factors, allowing for \ the fast
computation of a wide class of investment criteria. An example calibrated on
the real data shows that the optimal portfolio rules strongly depend on the
specification of the investment criterion. In particular, the rule is
sensitive to the inclusion of the correlation of portfolio value with
factors in the criterion. This supports the view that the correct choice of
the approximate investment criterion is crucial for optimal portfolio
allocation.

The remainder is organized as follows. Section \ref{Model} explains the
model. Section \ref{Result} formulates the main result. Section \ref%
{Application} specializes the result to a one-factor example, calibrates it
on the real data and presents numerical illustrations. And Section \ref%
{Conclusion} concludes. The proof of the main result is in Appendix.

\section{Model}

\label{Model}

\textbf{Assumption 1.} Securities follow a dynamic factor model.

\begin{align}
\frac{dS_i(t)}{S_i(t)}&=(a+AX(t))_i dt+\sum_{k=1}^{m+n}\sigma_{ik}dW_{k}(t), 
\hspace{1cm} i=1,2,\dots,m, \\
dX(t)&=BX(t)dt+\Lambda dW(t),
\end{align}

\noindent where $W(t)$ is a $\Re^{m+n}$ valued standard Brownian motion
process, $X(t)$ is $\Re^n$ valued factor process. It is assumed that $n
\times n$ matrix $B$ is stable. This coincides with the model of \cite%
{bielecki_pliska99} and \cite{bielecki_pliska_sherris00} with the exception
that the factor process $X(t)$ is normalized to have zero mean.

The process for investor's wealth $U$ is 
\begin{equation}
dU=\frac{IU}{S}dS,
\end{equation}%
where $I$ is the investment vector measured in shares of wealth.

\textbf{Assumption 2.} $I$ is a linear function of factors. 
\begin{equation}
I=h+HX,
\end{equation}%
where $H$ is a constant $m\times n$ matrix and $h$ is a constant $m\times 1$
vector.\footnote{%
This is not a very stringent assumption because non-linear functions of the
factors can always be added to the set of all factors to incorporate
non-linear dependence of investments on the factors.}

Thus, the investor's wealth $U$ follows the process 
\begin{equation}
dU=U(h+HX)^{\prime }\frac{dS}{S}
\end{equation}%
\noindent

\textbf{Assumption 3.} The investment criterion is 
\begin{equation}
W(h,H)=\liminf_{t\rightarrow \infty }\frac{1}{t}\{Eu(t)-\frac{\theta }{4}%
var(u(t))\}+\Gamma \lim_{t\rightarrow \infty }E(u(t)X(t)),  \label{target}
\end{equation}%
where $u(t)=:\ln U(t)-\ln U(0)$.

Note that this criterion includes the covariance of the logarithm of
portfolio value with the factors. An interpretation of this criterion is
that the investor cares not only about growth of the portfolio value and its
volatility but also about the ability of the portfolio to generate good
returns in bad times when his other income sources bring low returns.

\section{Result}

\label{Result} \noindent The following theorem shows how to compute
asymptotic moments that enter the investment criterion.

\begin{Theo}
\label{robport:theorem1} If assumptions 1 and 2 hold then 
\begin{align}
E(u)& =t[h^{\prime }a-\frac{1}{2}h^{\prime }\Sigma \Sigma ^{\prime }h+%
\mathrm{tr}(\Delta (H^{\prime }A-\frac{1}{2}H^{\prime }\Sigma \Sigma
^{\prime }H))], \\
\mathrm{Var}(u)& \sim t[YY^{\prime }+\mathrm{tr}(2SH^{\prime }A+(\Delta
-S)H^{\prime }\Sigma \Sigma ^{\prime }H)]+\mathrm{const}\text{\textrm{\ }}%
\mathrm{as}\text{ }t\rightarrow \infty , \\
E(uX)& \sim B^{-1}[\Delta (H^{\prime }\Sigma \Sigma ^{\prime }h-A^{\prime
}h-H^{\prime }a)-\Lambda \Sigma ^{\prime }h]\text{ }\mathrm{as}\text{ }%
t\rightarrow \infty .
\end{align}
\end{Theo}

Here $\Delta _{ij}=:E(x_{i}x_{j})$, vector $Y$ is defined as follows:

\begin{eqnarray}
Y=:(h^{\prime}\Sigma\Sigma^{\prime}H -h^{\prime}A -a^{\prime}H)B^{-1}\Lambda
+ h^{\prime}\Sigma,
\end{eqnarray}
and matrix $S$ is the symmetric solution of the equation 
\begin{eqnarray}
BS+SB^{\prime}=-2\Delta H^{\prime}A \Delta + \Delta
H^{\prime}\Sigma\Sigma^{\prime}H \Delta - 2\Lambda\Sigma^{\prime}H\Delta.
\end{eqnarray}

The proof of Theorem \ref{robport:theorem1} is in the Appendix. Let me
explain here its main idea for the case with only one factor $x$. Suppose
that we want to compute $\lim_{t\rightarrow\infty}E(u^i x^j)$ and we already
know $\lim_{t\rightarrow\infty}E(u^k x^l)$ for $k<i$ and for $k=i, l<j$.
Using Ito Lemma, we can write 
\begin{equation}
\begin{split}
d(u^i x^j)=&iu^{i-1}x^jdu+jx^{j-1}u^i dx \\
&+\frac{1}{2}\{i(i-1)u^{i-2}x^i(du)^2+j(j-1)u^ix^{j-2}(dx)^2+i j
u^{i-1}x^{j-1}dudx\}.
\end{split}%
\end{equation}
\noindent After the expressions for $du$ and $dx$ are substituted and the
expectations are taken, this equation is reduced to the following
differential equation on $E(u^i x^j)$: 
\begin{equation}
\frac{d}{dt}E(u^i x^j)=jBE(u^i x^j)+\dots,
\end{equation}
where $\dots$ includes only terms with already known asymptotic limits.
Since $B$ is stable (in this one factor example it means $B<0$), the
solution of this differential equation has an asymptotic limit that can be
easily calculated. Theorem 1 extends this method to the multi-factor case.



\section{Application}

\label{Application}

This section applies the general result to the particular example with the
parameters estimated from real data. First, Theorem \ref{robport:theorem1}
is specialized to the case with only one asset and one factor. Thus, $a$, $A$%
, $B$, $h$, and $H$ are scalars, $\Sigma=(\sigma ,\eta)$, and $\Lambda=(0
,\lambda)$.

In this setting, Theorem \ref{robport:theorem1} is reduced to the following

\begin{Coro}
\begin{align}
E(u)& =t\{ha-\frac{\lambda ^{2}}{2B}HA-\frac{\sigma ^{2}+\eta ^{2}}{2}(h^{2}-%
\frac{\lambda ^{2}}{2B}H^{2})\} \\
E(ux)& \sim \frac{\lambda ^{2}}{2B^{2}}\{hH(\eta ^{2}+\sigma ^{2})-Ah-Ha\}-%
\frac{\lambda \eta }{B}h \\
Var(u)& \sim t\{YY^{\prime }+2SHA+(\frac{\lambda ^{2}}{2B}-S)H^{2}(\sigma
^{2}+\eta ^{2})\}
\end{align}
\end{Coro}

where 
\begin{align*}
Y& =\{hH(\sigma ^{2}+\eta ^{2})-hA-aH\}\frac{1}{B}(0,\lambda )+h(\sigma
,\eta ) \\
S& =\frac{\lambda ^{2}}{4B^{2}}\{-2HA\frac{\lambda ^{2}}{2B}+\frac{\lambda
^{2}}{2B}H^{2}(\sigma ^{2}+\eta ^{2})-2H\lambda \eta \}.
\end{align*}

I have estimated the parameters of the stochastic process using the monthly
data from CRSP files. The effective stock return is taken as the difference
of the rate of return on S$\&$P composite index and the short-term interest
rate. The factor is the short-term interest rate, which is the rate on
3-month Treasury Bill. The results of estimation are summarized in the next
two tables.

\begin{minipage}[t]{6in}
Table 1. 
Process Estimates\footnote{\noindent Period is from January 1970 to December 2000. Rates and returns are measured in percentage terms; t-ratios are in parentheses}
\par

\begin{tabular}{|l|c|c|} \hline
&constant&r\\	\hline
dS/S&1.993&-1.177\\
&(3.505)&(-14.220)\\	\hline
dr&0.120&0.979\\
&(0.911)&(42.885)\\ \hline	
\end{tabular}
\end{minipage}

\begin{minipage}[t]{6in}
Table 2. 
Covariances of interest rate and stock return innovations\footnote{Period is from January 1970 to December 2000}

\begin{tabular}{|l|c|c|} \hline
&Stock Return&r\\	\hline
Stock Return&19.587&0.0553\\	\hline
r&0.0553&0.4006\\	\hline
\end{tabular}
\end{minipage}

\bigskip

Therefore, the parameters of the model can be calibrated as $a=0.01993$, $%
A=-0.01177$, $B=-0.021$, $\lambda =0.6329$, $\eta =0.000874$, and $\sigma
=0.044249$.

For the first illustration, assume that $h=1$, so that on average the
investor have all his funds in the stock. If the interest rate deviates from
its average value then the investor can move his assets to bonds. He can
also borrow additional funds and increase his position in the stock. The
first three graphs show how $Eu$, $Var(u)$, and $E(ux)$ change with respect
to changes in $H$.

\begin{figure}[tbp]
\caption{Dependence of $K=E(du)$ from $H$}
\label{robport:figure1}
\epsfbox{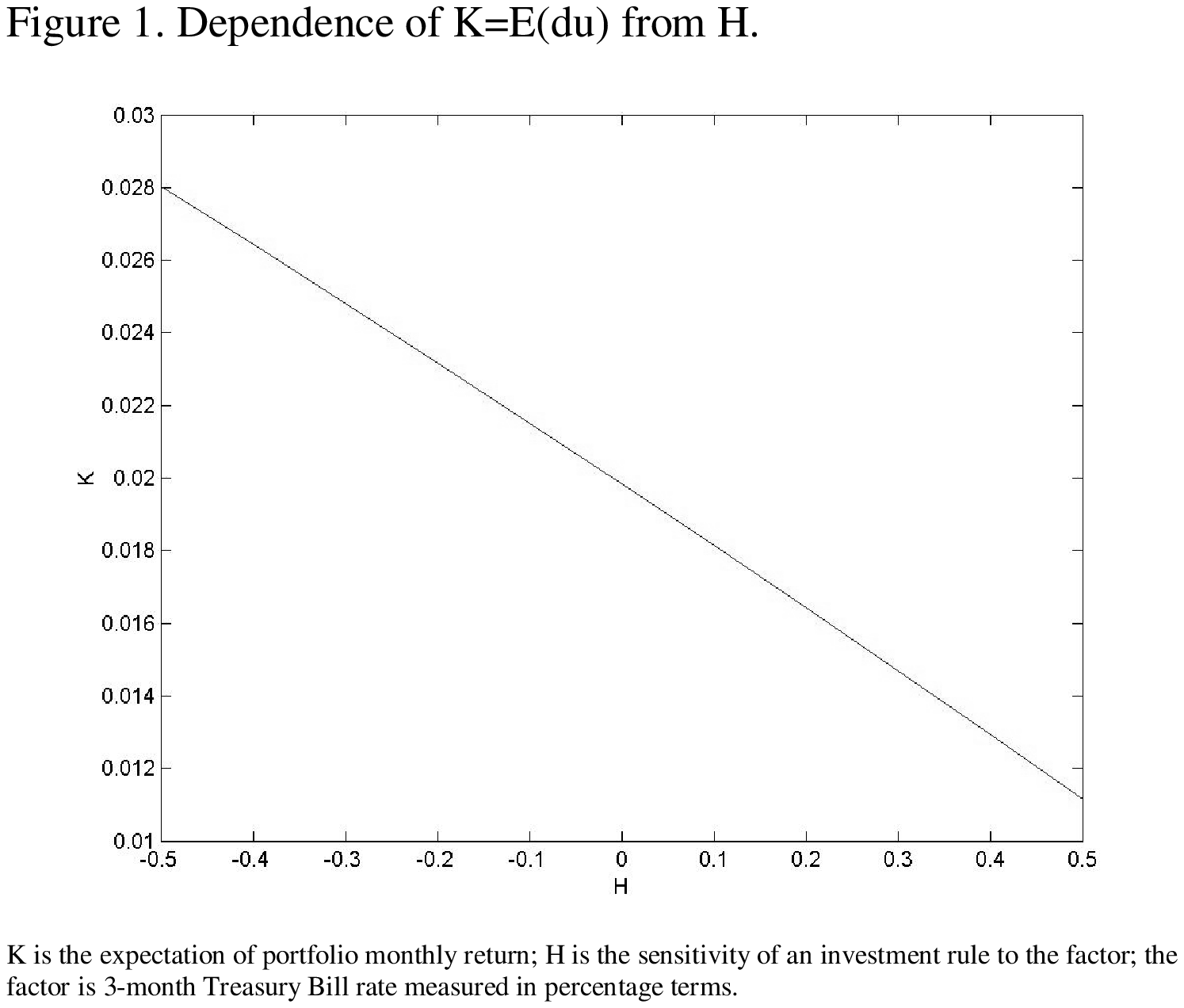}
\end{figure}

\begin{figure}[tbp]
\caption{Dependence of $P=\lim_{t\rightarrow\infty}E(Xu)$ from $H$}
\label{robport:figure2}
\epsfbox{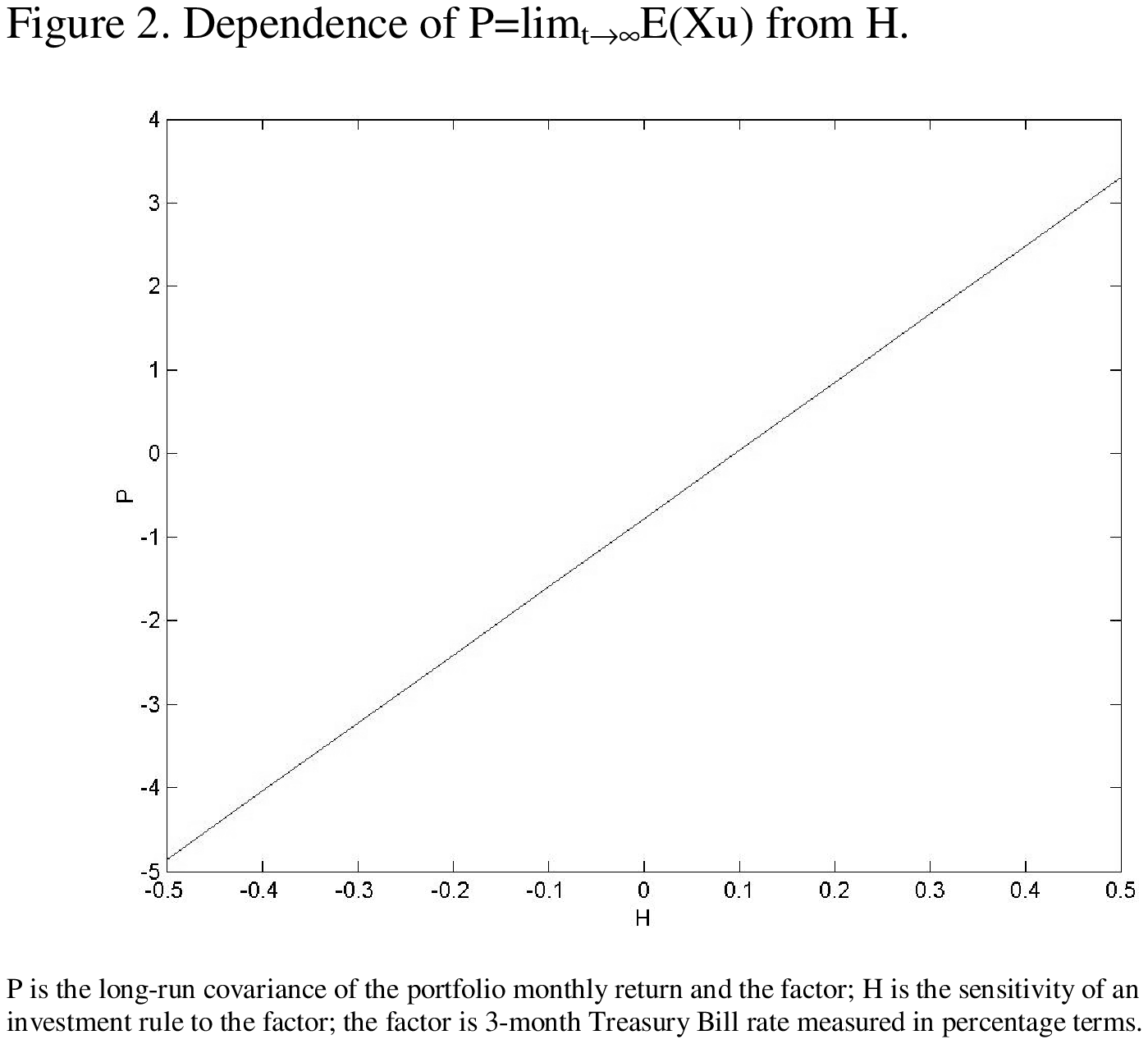}
\end{figure}

\begin{figure}[tbp]
\caption{Dependence of $R=\lim_{t\rightarrow\infty}(1/t)\Var(u)$ from $H$}
\label{robport:figure3}
\epsfbox{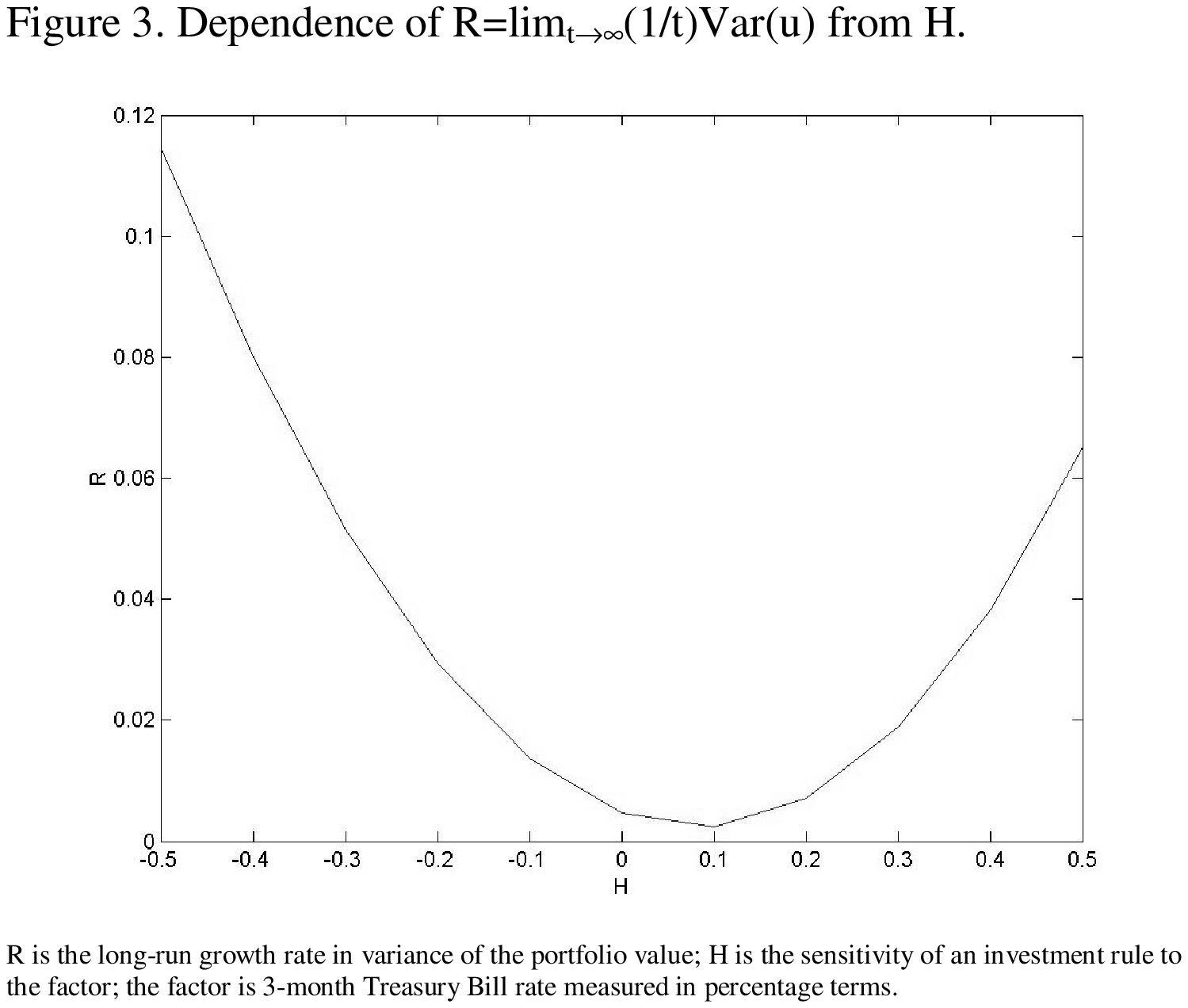}
\end{figure}


These graphs show that the strategy of investing more when the interest
rates are low increases the expected return of the portfolio. Expected
covariance of portfolio value with factors increases for large $H$.
Surprisingly, the volatility of the portfolio is minimized at a certain
non-zero level of $H$.

\begin{figure}[tbp]
\caption{Dependence of optimal $(h,H)$-strategy from risk-sensitivity
parameter $\protect\theta$}
\label{robport:figure4}
\epsfbox{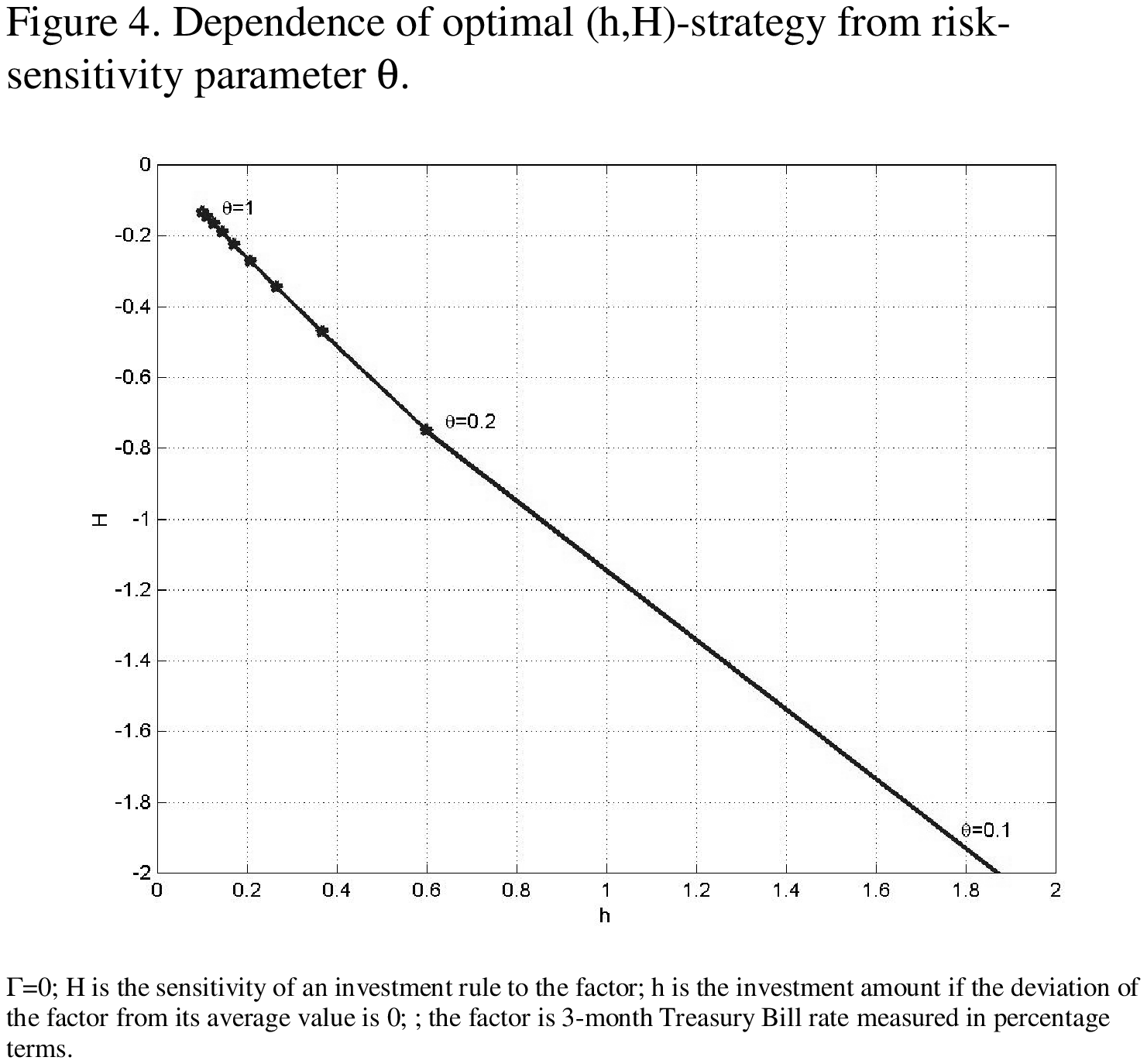}
\end{figure}

For the second illustration, assume that the target function (\ref{target})
is specified with factor-sensitivity $\Gamma=0$. Then the dependence of
optimal $(h,H)$ combination on the risk-sensitivity parameter $\theta$ is
illustrated on Figure \ref{robport:figure4}.

This picture shows that the increase in the risk-sensitivity leads to
smaller average investment in stocks and to smaller sensitivity of the
investment to changes in the interest rates. It is interesting to note,
however, that the ratio of $H$ to $h$, which can be interpreted as the
relative sensitivity of investments to interest rates is almost constant for
sufficiently large values of risk-sensitivity.

\begin{figure}[tbp]
\caption{Dependence of optimal $(h,H)$-strategy from risk-sensitivity
parameter $\protect\theta$ and factor-sensitivity parameter $\Gamma$}
\label{robport:figure5}
\epsfbox{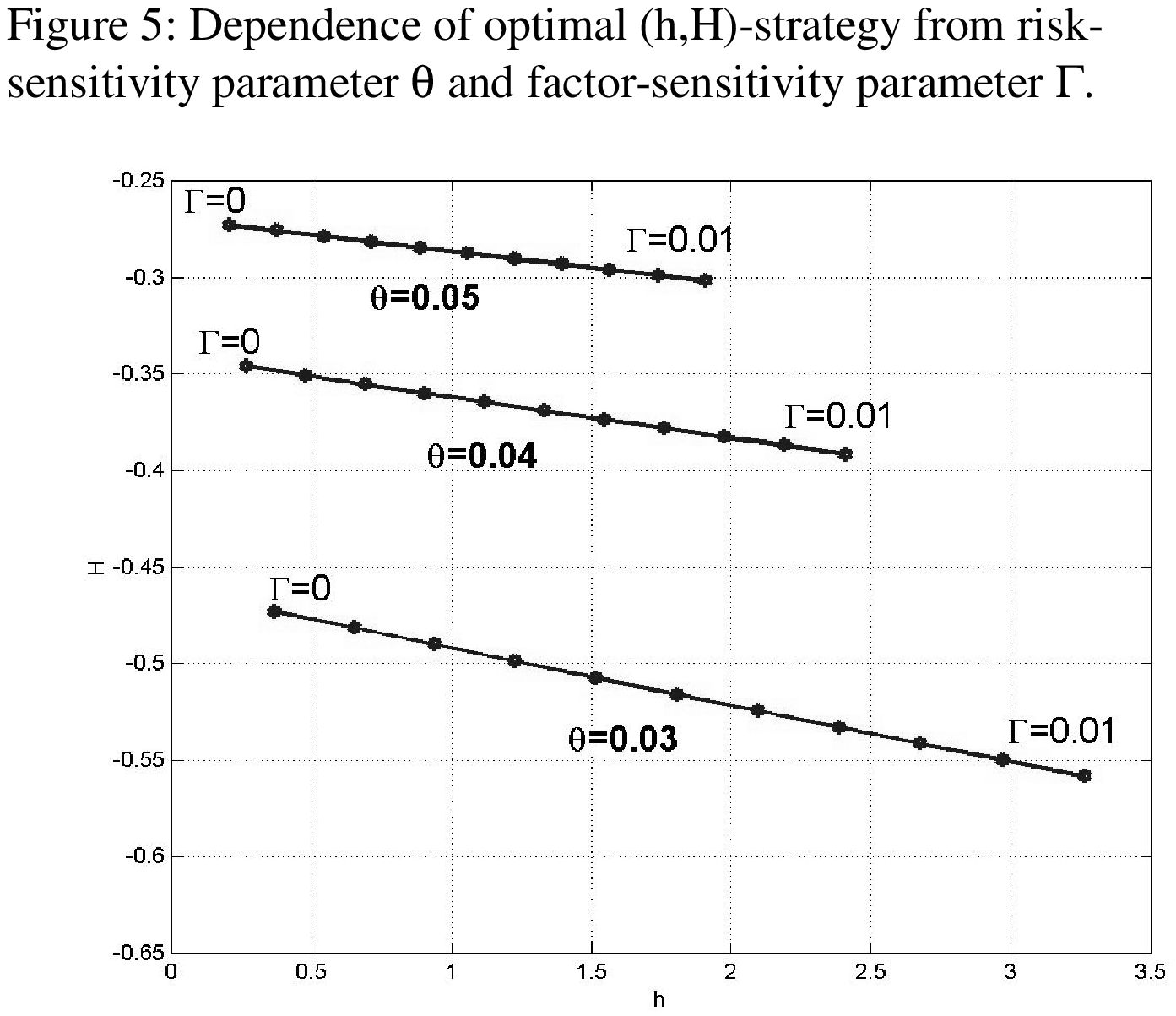}
\end{figure}

Now consider factor-sensitivity parameter $\Gamma $ changing from $0$ to $%
0.01$. The graph on Figure~\ref{robport:figure5} shows the optimal $(h,H)$
path. The greater sensitivity of the investment rule to the covariance of
the portfolio value with the factor leads to smaller average investment and
smaller sensitivity of the investment rule to the factor. It is interesting,
however, that the sensitivity of investment rule to factor is affected
smaller than the average investment amount.

\section{Conclusion}

\label{Conclusion}

This paper derives explicit formulas for asymptotic joint moments of the
logarithm of investor's wealth and factors, allowing for a fast maximization
of linear combinations of the moments. An application of this method to the
real data shows that the optimal strategy is very sensitive to the
factor-sensitivity parameter $\Gamma $ that measures dependence of the
investment criterion on the covariance of wealth with factors.

\appendix

\section{Proof of Theorem \ref{robport:theorem1}}

First, we are going to compute $Eu$: 
\begin{equation}
\begin{split}
du&=d(\ln U)=\frac{dU}{U}-\frac{1}{2}(\frac{dU}{U})^2 \\
&=(h+HX)^{\prime}((a+AX)dt+\Sigma dW) \\
&- \frac{1}{2}(h+HX)^{\prime}\Sigma\Sigma^{\prime}(h+HX)dt;
\end{split}%
\end{equation}
\begin{equation}  \label{Edu}
Edu=(h^{\prime}a+\text{tr}(A\Delta H^{\prime}))dt - \frac{1}{2}%
(h^{\prime}\Sigma\Sigma^{\prime}h+ \text{tr}(\Sigma^{\prime}H\Delta
H^{\prime}\Sigma))dt.
\end{equation}

Integrating (\ref{Edu}) over time gives 
\begin{equation}
Eu=Kt,
\end{equation}

\noindent where 
\begin{equation}
K=:h^{\prime}a+\text{tr}(A\Delta H^{\prime})-\frac{1}{2}(h^{\prime}\Sigma%
\Sigma^{\prime}h+ \text{tr}(\Sigma^{\prime}H\Delta H^{\prime}\Sigma)).
\end{equation}

\noindent Our second goal is to compute $E(uX)$. 
\begin{equation}
\begin{split}
d(uX)&=Xdu+udX+dudX \\
&= X[(h+HX)^{\prime}((a+AX)dt+\Sigma dW) - \frac{1}{2}(h+HX)^{\prime}\Sigma%
\Sigma^{\prime}(h+HX)dt] \\
&+u(BXdt+\Lambda dW) \\
&+[(h+HX)^{\prime}((a+AX)dt+\Sigma dW) - \frac{1}{2}(h+HX)^{\prime}\Sigma%
\Sigma^{\prime}(h+HX)dt](BXdt+\Lambda dW)
\end{split}%
\end{equation}

\noindent Using $E(X)=E(X^3)=0$, we get 
\begin{eqnarray}
Ed(uX)&=&dt(\Delta A^{\prime}h+ \Delta H^{\prime}a - \Delta H^{\prime}\Sigma
\Sigma^{\prime}h + BE(uX)+\Lambda \Sigma^{\prime}h)
\end{eqnarray}
Denoting $E(uX)$ as $P(t)$, we note that it satisfies the following
differential equation 
\begin{equation}
\frac{dP}{dt}=BP+\Delta A^{\prime}h+ \Delta H^{\prime}a - \Delta
H^{\prime}\Sigma \Sigma^{\prime}h + \Lambda \Sigma^{\prime}h.
\end{equation}
Since $B$ is stable, 
\begin{equation}
P(t)\rightarrow B^{-1}[\Delta(H^{\prime}\Sigma\Sigma^{\prime}h - A^{\prime}h
- H^{\prime}a)-\Lambda \Sigma^{\prime}h] \text{ as } t\rightarrow \infty.
\end{equation}
Thus, we computed $P=:\lim_{t\rightarrow \infty} P(t)$, which is the
asymptotic covariance of the portfolio and factors.

The next step is to compute $E(uXX^{\prime})$. 
\begin{equation}
\begin{split}
d(uXX^{\prime})&=duXX^{\prime}+u(dX)X^{\prime}+uXdX^{\prime}+dudX
X^{\prime}+ (du)XdX^{\prime}+udXdX^{\prime} \\
&=[(h+HX)^{\prime}((a+AX)dt+\Sigma dW) - \frac{1}{2}(h+HX)^{\prime}\Sigma%
\Sigma^{\prime}(h+HX)dt]XX^{\prime} \\
&+u(BXdt+\Lambda dW)X^{\prime} \\
&+uX(BXdt+\Lambda dW)^{\prime} \\
&+[(h+HX)^{\prime}((a+AX)dt+\Sigma dW) - \frac{1}{2}(h+HX)^{\prime}\Sigma%
\Sigma^{\prime}(h+HX)dt](BXdt+\Lambda dW)X^{\prime} \\
&+[(h+HX)^{\prime}((a+AX)dt+\Sigma dW) - \frac{1}{2}(h+HX)^{\prime}\Sigma%
\Sigma^{\prime}(h+HX)dt]X(BXdt+\Lambda dW)^{\prime} \\
&+u(BXdt+\Lambda dW)(BXdt+\Lambda dW)^{\prime}
\end{split}%
\end{equation}

\begin{equation}
\begin{split}
Ed(uXX^{\prime})&=dt[h^{\prime}a\Delta+E(\text{tr}(AXX^{\prime}H^{%
\prime})XX^{\prime}) \\
&-\frac{1}{2}h^{\prime}\Sigma\Sigma^{\prime}h \Delta - \frac{1}{2}E(\text{tr}%
(\Sigma^{\prime}HXX^{\prime}H^{\prime}\Sigma)XX^{\prime}) \\
&+BE(uXX^{\prime})+E(uXX^{\prime})B^{\prime}+2\Lambda\Sigma^{\prime}H
\Delta+ E(u)\Lambda \Lambda^{\prime}]
\end{split}%
\end{equation}

\noindent Asymptotically, $E(uXX^{\prime})\sim Rt+S$, and the coefficients $%
R,S$ can be found from the equation:

\begin{equation}
\begin{split}
R&=h^{\prime}a\Delta+E(\text{tr}(AXX^{\prime}H^{\prime})XX^{\prime})-\frac{1%
}{2}h^{\prime}\Sigma\Sigma^{\prime}h \Delta - \frac{1}{2}E(\text{tr}%
(\Sigma^{\prime}HXX^{\prime}H^{\prime}\Sigma)XX^{\prime}) \\
&+B(Rt+S)+(Rt+S)B^{\prime}+2\Lambda\Sigma^{\prime}H \Delta+ K\Lambda
\Lambda^{\prime}t,
\end{split}%
\end{equation}

\noindent which reduces to two equations on $R$ and $S$, 
\begin{align}  \label{R}
0&=BR + RB^{\prime}+ K\Lambda \Lambda^{\prime}, \\
\begin{split}  \label{S}
R&=h^{\prime}a\Delta+E(\text{tr}(AXX^{\prime}H^{\prime})XX^{\prime})-\frac{1%
}{2}h^{\prime}\Sigma\Sigma^{\prime}h \Delta - \frac{1}{2}E(\text{tr}%
(\Sigma^{\prime}HXX^{\prime}H^{\prime}\Sigma)XX^{\prime}) \\
& +BS+SB^{\prime}+2\Lambda\Sigma^{\prime}H \Delta.
\end{split}%
\end{align}

\noindent Since $\Delta$ satisfies $B\Delta + \Delta B^{\prime}+ \Lambda
\Lambda^{\prime}= 0$, (\ref{R}) implies 
\begin{equation}
R=K\Delta
\end{equation}
and, therefore, (\ref{S}) can be rewritten as 
\begin{equation}
\begin{split}  \label{S1}
K\Delta&=h^{\prime}a\Delta+E(\text{tr}(AXX^{\prime}H^{\prime})XX^{\prime})-%
\frac{1}{2}h^{\prime}\Sigma\Sigma^{\prime}h \Delta - \frac{1}{2}E(\text{tr}%
(\Sigma^{\prime}HXX^{\prime}H^{\prime}\Sigma)XX^{\prime}) \\
&+BS+SB^{\prime}+2\Lambda\Sigma^{\prime}H \Delta.
\end{split}%
\end{equation}

\noindent Using the identity 
\begin{equation}
E(x_i x_j x_k
x_l)=\Delta_{ij}\Delta_{kl}+\Delta_{ik}\Delta_{jl}+\Delta_{il}\Delta_{jk},
\end{equation}
we can further reduce (\ref{S1}) to 
\begin{equation}
BS+SB^{\prime}=-2\Delta H^{\prime}A \Delta + \Delta
H^{\prime}\Sigma\Sigma^{\prime}H \Delta - 2\Lambda\Sigma^{\prime}H\Delta.
\end{equation}

For later use, we also need to compute $E(uX^{\prime}QX)$ where $Q$ is an
arbitrary matrix. Because of the identity 
\begin{equation}
E(uX^{\prime}QX)=\text{tr}(E(uXX^{\prime})Q),
\end{equation}
we have 
\begin{equation}
E(uX^{\prime}QX)\sim \text{tr}(RQ)t+\text{tr}(SQ) \text{ when }
t\rightarrow\infty
\end{equation}

The next step is to compute $E(u^2)$. 
\begin{equation}
\begin{split}
d(u^2)&=2u(du)+ (du)^2 \\
&=2u[(h+HX)^{\prime}((a+AX)dt+\Sigma dW) - \frac{1}{2}(h+HX)^{\prime}\Sigma%
\Sigma^{\prime}(h+HX)dt] \\
&+(h+HX)^{\prime}\Sigma\Sigma^{\prime}(h+HX)dt
\end{split}%
\end{equation}

\begin{equation}
\begin{split}
Ed(u^2)&=2dt[h^{\prime}a Eu+E(uX^{\prime})
H^{\prime}a+h^{\prime}AE(uX)+E(uX^{\prime}H^{\prime}AX) \\
&-\frac{1}{2}(h^{\prime}\Sigma\Sigma^{\prime}h E(u)+h^{\prime}\Sigma
\Sigma^{\prime}HE(uX)+E(uX^{\prime})H^{\prime}\Sigma\Sigma^{\prime}h+E(uX^{%
\prime}H^{\prime}\Sigma \Sigma^{\prime}H X)) \\
&+E((h+HX)^{\prime}\Sigma\Sigma^{\prime}(h+HX))] \\
&=2dt[h^{\prime}a K t + P^{\prime}H^{\prime}a + h^{\prime}AP+\text{tr}%
(RH^{\prime}A)t+\text{tr}(SH^{\prime}A) \\
&-\frac{1}{2}(h^{\prime}\Sigma\Sigma^{\prime}h Kt +h^{\prime}\Sigma
\Sigma^{\prime}HP + P^{\prime}H^{\prime}\Sigma \Sigma^{\prime}h + \text{tr}%
(RH^{\prime}\Sigma \Sigma^{\prime}H) t + \text{tr}(SH^{\prime}\Sigma
\Sigma^{\prime}H)) \\
&+\frac{1}{2}E((h+HX)^{\prime}\Sigma\Sigma^{\prime}(h+HX))]
\end{split}%
\end{equation}

From this we have, 
\begin{equation}
\begin{split}
Eu^2&=\text{const}+ t^2[h^{\prime}a K + \text{tr}(RH^{\prime}A) -\frac{1}{2}%
(h^{\prime}\Sigma\Sigma^{\prime}h K + \text{tr}(RH^{\prime}\Sigma
\Sigma^{\prime}H) )] \\
&+t[2(P^{\prime}H^{\prime}a + h^{\prime}AP + \text{tr}(SH^{\prime}A))
-(h^{\prime}\Sigma \Sigma^{\prime}HP + P^{\prime}H^{\prime}\Sigma
\Sigma^{\prime}h +\text{tr}(SH^{\prime}\Sigma \Sigma^{\prime}H)) \\
&+E((h+HX)^{\prime}\Sigma\Sigma^{\prime}(h+HX))].
\end{split}%
\end{equation}

\noindent Finally, 
\begin{equation}
\begin{split}
Var(u)&=E(u^2)-E(u)^2 \\
&=\text{const}+2(P^{\prime}H^{\prime}a + h^{\prime}AP + tr(SH^{\prime}A))t \\
&-(h^{\prime}\Sigma \Sigma^{\prime}HP + P^{\prime}H^{\prime}\Sigma
\Sigma^{\prime}h +\text{tr}(SH^{\prime}\Sigma \Sigma^{\prime}H))t \\
&+(h^{\prime}\Sigma\Sigma ^{\prime}h+\text{tr}(\Sigma^{\prime}H\Delta
H^{\prime}\Sigma)t
\end{split}%
\end{equation}

Using the definition of $P$ this can be manipulated into 
\begin{equation}
\begin{split}
Var(u)& =t[((h^{\prime }\Sigma \Sigma ^{\prime }H-h^{\prime }A-a^{\prime
}H)B^{-1}\Lambda +h^{\prime }\Sigma )((h^{\prime }\Sigma \Sigma ^{\prime
}H-h^{\prime }A-a^{\prime }H)B^{-1}\Lambda +h^{\prime }\Sigma )^{\prime } \\
& +\text{tr}(2SH^{\prime }A+(\Delta -S)H^{\prime }\Sigma \Sigma ^{\prime
}H)]+\text{const}
\end{split}%
\end{equation}

\newpage

\bibliographystyle{plain}
\bibliography{comtest}

\end{document}